\documentclass[11pt]{article}

\usepackage{parskip}
\usepackage[T1]{fontenc}
\newcommand{\blind}{0}

\usepackage[margin=1in]{geometry}
\usepackage{setspace}

\usepackage{amsmath,amssymb,amsthm,bm}
\usepackage{graphicx,psfrag,epsf}
\usepackage{epigraph,xcolor}
\usepackage{enumerate}
\usepackage{url}
\usepackage{amsmath,latexsym,amssymb,bm}
\usepackage{hyperref}
\usepackage{bbm}
\usepackage{mathtools}
\usepackage{comment}
\usepackage{mathrsfs}
\hypersetup{colorlinks=true,citecolor=blue,
	pdfpagemode=FullScreen,linktocpage=true}
\usepackage{cleveref}

\usepackage{natbib}
\newcommand{\RNum}[1]{\uppercase\expandafter{\romannumeral #1\relax}}

\newcommand{\paren}[1]{\mathopen{}\left( {#1}_{{}_{}}\,\negthickspace\right)\mathclose{}}

\numberwithin{equation}{section}

\newtheoremstyle{general}
{3mm} 
{3mm} 
{\it} 
{} 
{\bfseries} 
{.} 
{.5em} 
{} 

\theoremstyle{general}

\newtheorem{lemma}{Lemma}
\newtheorem{theorem}{Theorem}
\newtheorem{corollary}{Corollary}

\newtheorem{claim}{Claim}

\makeatletter
\renewenvironment{proof}[1][\proofname]{\par
    \pushQED{\qed}%
    \normalfont \topsep6\p@\@plus6\p@\relax
    \trivlist
    \item\relax{
        \bfseries
        #1\@addpunct{.}}\hspace\labelsep\ignorespaces
    }{%
     \popQED\endtrivlist\@endpefalse
     }
\makeatother

\begin{document}


\if0\blind
{
  \title{\bf Feedback Interacting Urn Models}
  \author{Krishanu Maulik \thanks{Email: krishanu@isical.ac.in}\hspace{.2cm} \\
  Indian Statistical Institute Kolkata \\ 203 Barrackpore Trunk Road, Kolkata, West Bengal 700108, India. \\
  and \\
  Manit Paul\thanks{Email: manit282000@gmail.com}\hspace{.2cm} \\
  Indian Statistical Institute Kolkata \\ 203 Barrackpore Trunk Road, Kolkata, West Bengal 700108, India.}
  \maketitle
} \fi

\if1\blind
{
  \bigskip
  \bigskip
  \bigskip
  \begin{center}
    {\LARGE\bf A Special Interacting Urn Model: \\ Some Interesting Consequences}
\end{center}
  \medskip
} \fi

\bigskip
\begin{abstract}
We introduce and discuss a special type of feedback interacting urn model with deterministic interaction. This is a generalisation of the very well known \cite{polyaurn} urn model. In our model, balls are added to a particular urn depending on the replacement matrix of that urn and the color of ball chosen from some other urn. This urn model can help in studying how various interacting models might behave in real life in the long run. We have also introduced a special type of interacting urn model with non-deterministic interaction and studied its behaviour. Furthermore, we have provided some nice examples to illustrate the various consequences of these interacting urn models.
\end{abstract}

\noindent%
\vfill





\newpage

\section{Introduction}
\label{sec:introduction}
We devise a toy model inspired by the following real life examples. Let us consider two interacting stock markets, A and B. Suppose initially the stock market-A is assumed to have bullish tendency while stock market-B is assumed to have bearish tendency. These two stock markets operate alternately in mutually exclusive hours, for example, Tokyo stock market and New York stock market. Now due to the interaction between these two stock markets, each market may get influenced by the other. Our model is motivated by a natural question that might arise here: whether it is plausible for both the stock markets to change their nature in the long run due to the influence of the other. Another related example is in the context of the interaction between political preferences of two neighbouring states in USA during their elections. We know the states in USA hold state governor elections every four years but not necessarily simultaneously. Election results of one state can influence that of neighbouring state holding elections in a different cycle. Thus we can study two neighbouring states running different election cycles with tilt towards Republicans and Democrats respectively and analyse whether their political leanings can go through a change.
We devise a toy model inspired by the following real life examples. Let us consider two interacting stock markets, A and B. Suppose initially the stock market-A is assumed to have bullish tendency while stock market-B is assumed to have bearish tendency. These two stock markets operate alternately in mutually exclusive hours, for example, Tokyo stock market and New York stock market. Now due to the interaction between these two stock markets, each market may get influenced by the other. Our model is motivated by a natural question that might arise here: whether it is plausible for both the stock markets to change their nature in the long run due to the influence of the other. Another related example is in the context of the interaction between political preferences of two neighbouring states in USA during their elections. We know the states in USA hold state governor elections every four years but not necessarily simultaneously. Election results of one state can influence that of neighbouring state holding elections in a different cycle. Thus we can study two neighbouring states running different election cycles with tilt towards Republicans and Democrats respectively and analyse whether their political leanings can go through a change.

We shall model these phenomena using urn models. Urn model is one of the simplest and most useful models considered in probability theory. There have been various extensions and generalizations of this model since its introduction by \cite{polyaurn}. In this work we shall consider a special extension of this model with multiple urns involving feedback interaction. We have considered cases where the interaction is deterministic and also cases where it is not. Similar interacting urn models have previously been studied in \cite{chen2014new}, \cite{launay2012generalized}, \cite{siegmund2005urn}, \cite{kaur2019interacting}.

To begin with, we consider a 2-urn-2-color model. In this model, we are given two urns, both of which contain balls of two colors, say red and black. We shall motivate our model by the stock market example. However a similar explanation will hold for the political behaviour. The color red will correspond to the bearish tendency and the color black will be for the bullish tendency. The two stock markets may interact with each other in various ways. For example a trader in the market-B (which is initially bearish) will observe the outcome of the market-A in the previous session to decide his move. If the market-A has moved bullishly, the trader will put some weight for bullish nature, denoted by $\beta$, while deciding his trade. However, if an inherently bullish market-A behaves bearishly in the previous session, the trader in market-B (which is bearish itself) will be sure to make a bearish trade. A trade in the stock market-A under the influence of the stock market-B will be defined analogously. Thus we define two replacement matrices, $\boldsymbol{R_1},\boldsymbol{R_2}$, as follows.
\begin{equation}
\nonumber
\boldsymbol{R_1}=\begin{pmatrix} \alpha & 1-\alpha \\ 0 & 1 \end{pmatrix}, \;\;\;\;\; \boldsymbol{R_2}=\begin{pmatrix} 1 & 0 \\ 1-\beta & \beta  \end{pmatrix},
\end{equation}
where $0< \alpha,\beta<1$. Initially in the first urn, black is the dominant color and in the second urn, red is the dominant color. We start with one ball of dominant color in each of the urns to signify their initial behaviour. An iterative process motivated by the stock market and the electoral example is then performed on these two urns. A ball is chosen from the first urn and its color is noted. The color of the chosen ball will correspond to the type of move in the stock market-A. Based on the color being red (bearish move) or black (bullish move), we reinforce the second urn (corresponding to the stock market-B) according to $R_2$. It can be seen that corresponding to the red color (bearish move) we reinforce the red color (bearish move) only in the second urn. On the other hand corresponding to the black color (bullish move) both the colors are reinforced in the second urn with suitable weighing. Thereafter a ball is chosen from the second urn and a similar process is done. This process is then repeated, alternately applying it on both of the urns. Here $\alpha$ and $\beta$ will denote the amount of influence of the other market on the markets A and B respectively.

In this work, we say that a particular urn \textit{flips} if the proportion of the dominant color of the urn becomes less than that of the other color. For instance in this model, we say that the first urn flips in the limiting case if the proportion of red balls becomes more than that of black balls. The flipping of the second urn can be understood analogously. It shall be shown that, both the urns can not flip simultaneously in the limiting case. In particular we shall show in \Cref{subsec:2urn} that,
\begin{align*}
\text{first urn flips iff } & (\alpha,\beta)\in S_1=\{(\alpha,\beta)|0<\alpha,\beta<1 \And 2\alpha > 1+\alpha\beta\} \\
\intertext{and} \text{second urn flips iff } & (\alpha,\beta)\in S_2=\{(\alpha,\beta)|0<\alpha,\beta<1 \And2\beta > 1+\alpha\beta\}.
\end{align*}
In this paper we shall extend the model to $n$ urns containing balls of $n$ colors ($n\geq 2$) and different feedback interaction mechanisms shall be appropriately defined and analysed.

The paper is divided into the following sections. \Cref{sec:deterministic cyclic} introduces a feedback interacting urn model where the interaction is deterministic. We also compute the almost sure limits of the proportion vectors of the urns (i.e. the vector of the  proportion of each color present in an urn) in this section. \Cref{sec:Non-deterministic urn} considers a far more generalized interacting urn model where the interaction is non-deterministic. We further compute the almost sure limits of the proportion vectors of the urns for this generalized model. \Cref{sec:consequences} discusses some nice and interesting consequences of the theorems developed in the earlier sections.

\section{Feedback Interacting Urn Model with Deterministic Interaction}
\label{sec:deterministic cyclic}
Here, we consider a deterministic cyclic interacting urn model. In this model, there are $m$ urns. Each of the $m$ urns may contain balls of $N$ different colors. The replacement matrices of the $m$ urns are given by $\boldsymbol{R_1}, \boldsymbol{R_2}, ..., \boldsymbol{R_m}$, all of which are $N \times N$ matrices. Initially there is at least one ball in each of the $m$ urns. The following iterative process is then performed:
\begin{enumerate}
    \item A ball is chosen from the first urn. The color of the ball is noted. Then, balls are added to the second urn according to $\boldsymbol{R_2}$ based on the color of the ball chosen from the first urn.
    \item A ball is chosen from the second urn. The color of the ball is noted. Then, balls are added to the third urn according to $\boldsymbol{R_3}$ based on the color of the ball chosen from the second urn.
    \item $\cdots$
    \item A ball is chosen from the $m$-th urn. The color of the ball is noted. Then, balls are added to the first urn according to $\boldsymbol{R_1}$ based on the color of the ball chosen from the $m$-th urn.
    \item Go to step(1).
\end{enumerate}
Under this set-up, we have the following theorem.
\begin{theorem}\label{thm:deterministic interacting}
Suppose $\boldsymbol{R}$ is the following $Nm \times Nm$ matrix: \begin{equation}
    \boldsymbol{R}=\begin{pmatrix} 0&\boldsymbol{R_2}&0& \hdots &0 \\
    0&0&\boldsymbol{R_3}& \hdots &0\\
    \vdots & \vdots & \vdots & \ddots & \vdots\\

    0&0&0& \hdots &\boldsymbol{R_m}\\
    \boldsymbol{R_1}&0&0& \hdots &0
    \end{pmatrix}
\end{equation}
We make the following two assumptions: \textit{(i)} $\boldsymbol{R_i}$ is a balanced (each row sum is 1) replacement matrix with all the entries non-negative for i=1,2,$\cdots$,m and \textit{(ii)} $\boldsymbol{R}$ is an irreducible matrix. \\
Let $\boldsymbol{C_n^i}$ denote the composition vector of the $i$th urn at the $n$th stage, where $i=1,2,\cdots,m$. Define $\boldsymbol{A_k^i}=\frac{\boldsymbol{{C_{m(k-1)+i-1}^i}}}{k}$ for $i=1,2,\cdots,m$. If we have the set-up described as above, then, we have the following result:
\begin{equation}
\nonumber
    \boldsymbol{A_t^i} \xrightarrow{a.s.} \boldsymbol{\pi_{(R_{i+1}R_{i+2}\cdots R_{m}R_1\cdots R_i)}} \;\; \text{as}\;\; t\xrightarrow{}\infty ,
\end{equation}
where  $\boldsymbol{\pi_{(R_{i+1}R_{i+2}\cdots R_{m}R_1...R_i)}}$ is the left eigen-vector of the matrix $\boldsymbol{R_{i+1}R_{i+2}\cdots R_{m}R_1...R_i}$ corresponding to the eigen-value 1. Moreover all the entries of $\boldsymbol{\pi_{(R_{i+1}R_{i+2}\cdots R_{m}R_1...R_i)}}$ are strictly positive.
\end{theorem}

\begin{proof}[Proof of \Cref{thm:deterministic interacting}]
We shall require the following lemma taken from \cite{gouet1997strong} to prove this theorem.
\begin{lemma} \label{lem:gouet}
Let $\boldsymbol{B}=((b_{ij}))$ be an irreducible non-negative $p \times p$ matrix and $\displaystyle \sum_jb_{ij}=1\;\;for\;i=1,2,\cdots,p$. Let $\boldsymbol{u_n}=(\boldsymbol{u_{n1}},\cdots,\boldsymbol{u_{np}})$ be a sequence of non-negative normalized (sum of components=1) vectors such that $\displaystyle \boldsymbol{u_{n+1}}-\overline{\boldsymbol{u_n}}\boldsymbol{B} \xrightarrow{}0$, then if $\boldsymbol{u}$ is the dominant eigen-vector of $\boldsymbol{B}$, we have the following convergence:
\begin{equation}
    \label{3.20}
    \boldsymbol{u_n}\xrightarrow{}\boldsymbol{u}.  \nonumber
\end{equation}
\end{lemma}
 We can now proceed to prove this theorem. We can re-frame the given model in the following way. We combine all the $m$ urns into a single urn. The replacement matrix of the combined urn is the $Nm \times Nm$ matrix $\boldsymbol{R}$, as given in \Cref{thm:deterministic interacting}. In the combined urn, the balls of first $N$ colors correspond to the balls in Urn-1, the balls of the next $N$ colors correspond to the balls in Urn-2 and so on. Let $\boldsymbol{C_n^i}$ denote the composition of the $i$th urn in the $n$th stage. Also, let $\displaystyle \boldsymbol{C_t}=(\boldsymbol{C_t^1},\boldsymbol{C_t^2},\cdots,\boldsymbol{C_t^m})$. In other words, $\boldsymbol{C_t}$ is the composition vector of the combined urn at the $t$-th stage. Let us denote by $\boldsymbol{\chi_n}$, the vector to indicate which color is chosen at the $n$th stage from the combined urn. Note that, both $\boldsymbol{C_n}$ and $\boldsymbol{\chi_n}$ are $Nm$ dimensional row vectors. Hence, $\boldsymbol{\chi_n}=\boldsymbol{e_i}$ if the $i$-th color is chosen at the $n$-th stage, where $i=1,2,3,\cdots,Nm$ ($\{\boldsymbol{e_1}, \boldsymbol{e_2}, \boldsymbol{e_3},\cdots,\boldsymbol{e_{Nm}}\}$ is the canonical basis of $\mathbb{R}^{Nm}$). We need to keep in mind that in the first step, the ball is chosen from urn-1, in the second step, the ball is chosen from urn-2 and so on in a cyclic manner. Suppose, $\mathscr{F}_n$ is the sigma field generated by the set $\{\boldsymbol{\chi_i} \; for \; i=1,2,\cdots,n-1\}$. Let us describe the distribution of $\boldsymbol{\chi_n}$. For $n=mk$ where $k \in \mathbb{Z}^+$, we have,
\begin{equation}
    \label{eq:distribution of chi, deterministic, m}
    \mathbb{P}(\boldsymbol{\chi_n}=\boldsymbol{e_i}|\mathscr{F}_n)=\begin{cases} & \frac{\boldsymbol{C_{n-1,i}}}{k+1} \; \text{for} \; i=(m-1)N+1,(m-1)N+2,\cdots,mN, \\ & 0 \;\; \text{otherwise} .\end{cases}
\end{equation}
Similarly, when $n=mk+r$ where $r=1,2,...,m-1$ and $k \in \mathbb{Z}^+ \cup \{0\} $, we have:
\begin{equation}
    \label{eq:distribution of chi, deterministic}
     \mathbb{P}(\boldsymbol{\chi_n}=\boldsymbol{e_i}|\mathscr{F}_n)=\begin{cases} & \frac{\boldsymbol{C_{n-1,i}}}{k+1} \; \text{for} \; i=(r-1)N+1,(r-1)N+2,\cdots,rN, \\ & 0 \;\; \text{otherwise} .\end{cases}
\end{equation}
We shall now apply the Martingale version of Second Borel-Cantelli Lemma separately $m$ times on $\boldsymbol{\chi_{n,i}}$ depending on whether $n$ is equal to $mk+1, mk+2,\cdots, mk$ for some $k \in \mathbb{Z}^+ \cup \{0\}$. (Note that if $\boldsymbol{V_n}$ is a vector, then $\boldsymbol{V_{n,i}}$ denotes the $i$-th component of $\boldsymbol{V_n}$). We get the following $m$ equations using the distribution of $\boldsymbol{\chi_n}$ given in \eqref{eq:distribution of chi, deterministic, m} and \eqref{eq:distribution of chi, deterministic}. For $n=mk+r$ where $r=1,2,\cdots,m$ , we have,
\begin{equation}
    \frac{\sum_{k=0}^t \boldsymbol{\chi_{n,i}}}{\sum_{k=0}^t {\frac{\boldsymbol{C_{n-1,i}}}{k+1}}}\xrightarrow{a.s.} 1 \; \text{on} \; \left\{ \sum_{k=0}^\infty {\frac{\boldsymbol{C_{n-1,i}}}{k+1}}=\infty \right\},
\end{equation}
for $i=(r-1)N+1,(r-1)N+2,\cdots,rN$. We have, $\displaystyle \sum_{k=0}^\infty {\boldsymbol{C_{mk+r,i}}}/(k+1) \geq \sum_{k=0}^\infty {1}/(k+1)=\infty \; for \; r=1,2,\cdots,m$. Thus, we get the following $m$ equations. For $n=mk+r$ where $r=1,2,...,m$ , we have,
\begin{equation}
     \label{eq:bc on distrn chi}
     \sum_{k=0}^t\boldsymbol{\chi_{mk+r,i}}=\paren{1+o_i\paren{t}}\sum_{k=0}^t\frac{\boldsymbol{C_{mk+r-1,i}}}{k+1}.
\end{equation}
Now, let us denote by $\boldsymbol{C_n^1}$, the sub-vector of $\boldsymbol{C_n}$ formed by its first $N$ entries; by $\boldsymbol{C_n^2}$, the sub-vector of $\boldsymbol{C_n}$ formed by its next $N$ entries and so on. Similarly, we define $\boldsymbol{\chi_n^1} , \boldsymbol{\chi_n^2},\cdots,\boldsymbol{\chi_n^m}$. Thus, we can rewrite the $m$ equations in \eqref{eq:bc on distrn chi} in the following manner. For $n=mk+r$ where $r=1,2,\cdots,m$ , we have,
\begin{equation}
    \sum_{k=0}^t \boldsymbol{\chi_{mk+r}^1}=\sum_{k=0}^t \frac{\boldsymbol{C_{mk+r-1}^1}}{k+1}.\paren{\boldsymbol{I_N} + \boldsymbol{o^r\paren{t}}},
\end{equation}
where, $\boldsymbol{o^r(t)}$ is the $\displaystyle N \times N$ diagonal matrix with the diagonal entries, $o_{(r-1)N+1}(t),o_{(r-1)N+2}(t),\cdots,o_{rN}(t)$. We know that for $n\geq1$ the urn composition evolves as :
\begin{equation}
\label{eq:evolution of urn composition}
    \boldsymbol{C_{n}}=\boldsymbol{C_0}+\paren{\sum_{j=1}^n\boldsymbol{\chi_j}}\boldsymbol{R} .
\end{equation}
There are $m$ cases depending on whether $n$ is equal to $mk+1, mk+2,\cdots, mk+m$ for some $k \in \mathbb{Z}^+ \cup \{0\}$. We shall show the detailed calculation for only $n=mk$ where $k \in \mathbb{Z}^+$. The calculations for the other cases can be done quite similarly. We have,
\begin{equation}
    \label{eq:case m}
    \begin{split}
    \sum_{j=1}^n\boldsymbol{\chi_j} &=\paren{\sum_{k=0}^{t-1}\boldsymbol{\chi_{mk+1}^1},\sum_{k=0}^{t-1}\boldsymbol{\chi_{mk+2}^2},\hdots,\sum_{k=1}^t\boldsymbol{\chi_{mk}^m}}\\ &=\paren{\sum_{k=0}^{t-1}\frac{\boldsymbol{C_{mk}^1}}{k+1}\paren{\boldsymbol{I_N}+\boldsymbol{o^1\paren{n}}},\sum_{k=0}^{t-1}\frac{\boldsymbol{C_{mk+1}^2}}{k+1}\paren{\boldsymbol{I_N}+\boldsymbol{o^2\paren{n}}},\hdots,\sum_{k=1}^t\frac{\boldsymbol{C_{mk-1}^m}}{k+1}\paren{\boldsymbol{I_N}+\boldsymbol{o^m\paren{n}}}} \\ &=\paren{\sum_{k=0}^{t-1}\frac{\boldsymbol{C_{mk}^1}}{k+1},\sum_{k=0}^{t-1}\frac{\boldsymbol{C_{mk+1}^2}}{k+1},\hdots,\sum_{k=1}^t\frac{\boldsymbol{C_{mk-1}^m}}{k+1})\paren{\boldsymbol{I_{Nm}}+\boldsymbol{o\paren{n}}}},
    \end{split}
\end{equation}
where, $\boldsymbol{o(n)}$ is the $\displaystyle Nm \times Nm$ diagonal matrix with the diagonal entries, $\boldsymbol{o^1(n)},\boldsymbol{o^2(n)},\cdots,\boldsymbol{o^m(n)}$. From \eqref{eq:evolution of urn composition} and \eqref{eq:case m}, we have,
\begin{equation}
    \boldsymbol{C_{mt}}-\boldsymbol{C_0}=\paren{\boldsymbol{I_{Nm}}+\boldsymbol{o\paren{n}}}\paren{\sum_{k=0}^{t-1}\frac{\boldsymbol{C_{mk}^1}}{k+1},\sum_{k=0}^{t-1}\frac{\boldsymbol{C_{mk+1}^2}}{k+1},\cdots,\sum_{k=1}^t\frac{\boldsymbol{C_{mk-1}^m}}{k+1}}\boldsymbol{R}.
\end{equation}
Matching the first components of the vectors on the right and left side, we obtain,
\begin{equation}
\label{eq:case m evolution}
    \frac{\boldsymbol{C_{mt}^1}-\boldsymbol{C_0^1}}{t+1}=\paren{\boldsymbol{I_{N}}+\boldsymbol{o^1\paren{n}}}\frac{1}{t+1}\sum_{k=1}^t\frac{\boldsymbol{C_{mk-1}^m}}{k+1}\boldsymbol{R_1}.
\end{equation}
Now, we define, $\displaystyle \boldsymbol{A_k^i}={\boldsymbol{C_{m(k-1)+i-1}^i}}/{k}$ for $i=1,2,\cdots,m$. Note that, $\boldsymbol{A_k^i}$ is the proportion vector of the $i$th urn at the $k$th stage(i.e. the first entry in $\boldsymbol{A_k^i}$ denotes the fraction of balls of first color in urn-$i$ at the $k$th stage, the second entry denotes the fraction of balls of second color in urn-$i$ at the $k$th stage and so on). \\ Thus, the cesaro means will be, $\displaystyle \Bar{\boldsymbol{A_t^i}}=({\sum_{k=1}^t\boldsymbol{A_k^i}})/{t}$ for $i=1,2,\cdots,m$. Using these notations, equation \eqref{eq:case m evolution} boils down to :
\begin{equation}
    \boldsymbol{A_{t+1}^1}-\overline{\boldsymbol{A_t^m}}\boldsymbol{R_1} \xrightarrow{a.s.}0 \;\; \text{as} \; t\xrightarrow{} \infty.
\end{equation}
Similarly, we get $m$ convergence equations:
\begin{equation}
   \begin{split}
       &\boldsymbol{A_{t}^2}-\overline{\boldsymbol{A_{t}^1}}\boldsymbol{R_2} \xrightarrow{a.s.}0 \;\; \text{as} \; t\xrightarrow{} \infty, \\
    &\boldsymbol{A_{t}^3}-\overline{\boldsymbol{A_{t}^2}}\boldsymbol{R_3} \xrightarrow{a.s.}0 \;\; \text{as} \; t\xrightarrow{} \infty,  \\
    & \hdots \\
    &\boldsymbol{A_{t}^m}-\overline{\boldsymbol{A_{t}^{m-1}}}\boldsymbol{R_m} \xrightarrow{a.s.}0 \;\; \text{as} \; t\xrightarrow{} \infty, \\
    &\boldsymbol{A_{t+1}^1}-\overline{\boldsymbol{A_{t}^m}}\boldsymbol{R_1} \xrightarrow{a.s.}0 \;\; \text{as} \; t\xrightarrow{} \infty.
   \end{split}
\end{equation}
We can rewrite these $m$ equations as follows:
\begin{equation}
    \begin{split}
        &\boldsymbol{A_{t+1}^1}-\overline{\boldsymbol{A_{t}^m}}\boldsymbol{R_1} \xrightarrow{a.s.}0 \;\; \text{as} \; t\xrightarrow{} \infty, \\
&\boldsymbol{A_{t+1}^2}-\frac{\boldsymbol{A_{t+1}^1}}{t+1}\boldsymbol{R_2}-\frac{t}{t+1}\overline{\boldsymbol{A_{t}^1}}\boldsymbol{R_2} \xrightarrow{a.s.}0 \;\; \text{as} \; t\xrightarrow{} \infty, \\
&\boldsymbol{A_{t+1}^3}-\frac{\boldsymbol{A_{t+1}^2}}{t+1}\boldsymbol{R_3}-\frac{t}{t+1}\overline{\boldsymbol{A_{t}^2}}\boldsymbol{R_3} \xrightarrow{a.s.}0 \;\; \text{as} \; t\xrightarrow{} \infty, \\
&\hdots\\
&\boldsymbol{A_{t+1}^m}-\frac{\boldsymbol{A_{t+1}^{m-1}}}{t+1}\boldsymbol{R_m}-\frac{t}{t+1}\overline{\boldsymbol{A_{t}^{m-1}}}\boldsymbol{R_m} \xrightarrow{a.s.}0 \;\; \text{as} \; t\xrightarrow{} \infty .
    \end{split}
\end{equation}
Now, we shall combine all the $m$ equations into a single equation. Let $\displaystyle \boldsymbol{D_t}=(\boldsymbol{A_t^1},\boldsymbol{A_t^2},\cdots,\boldsymbol{A_t^m})$. So, on combining all the $m$ equations into a single equation, we get:
\begin{equation}
    \boldsymbol{D_{t+1}}\boldsymbol{G_t}-\overline{\boldsymbol{D_t}}\boldsymbol{H_t} \xrightarrow{a.s.}0\;\;\text{as}\;\;t\xrightarrow{}\infty,
\end{equation}
which implies that,
\begin{equation}
    \label{eq:combined equation}
     \boldsymbol{D_{t+1}}-\overline{\boldsymbol{D_t}}\boldsymbol{H_t}\boldsymbol{G_t}^{-1}\xrightarrow{a.s.}0\;\;\text{as}\;\;t\xrightarrow{}\infty,
\end{equation}
where,
\begin{equation}
    \label{eq:Gt,Ht}
       \boldsymbol{G_t}=\begin{pmatrix}\boldsymbol{I_N}&-\frac{\boldsymbol{R_2}}{t+1}&0&\hdots&0\\
        0&\boldsymbol{I_N}&-\frac{\boldsymbol{R_3}}{t+1}&\hdots&0\\
        \vdots&\vdots&\vdots&\ddots&\vdots\\
        0&0&0&\hdots&-\frac{\boldsymbol{R_m}}{t+1}\\
        0&0&0&\hdots&\boldsymbol{I_N}
        \end{pmatrix},\;
    H_t=\begin{pmatrix}0&\frac{t\boldsymbol{R_2}}{t+1}&0&\hdots&0\\
        0&0&\frac{t\boldsymbol{R_3}}{t+1}&\hdots&0\\
        \vdots&\vdots&\vdots&\ddots&\vdots\\
        0&0&0&\hdots&\frac{t\boldsymbol{R_m}}{t+1}\\
        \boldsymbol{R_1}&0&0&\hdots&0
        \end{pmatrix}.
\end{equation}
We observe that, $\displaystyle \boldsymbol{G_t}\xrightarrow{a.s.}\boldsymbol{I_{Nm}}$ as $t\xrightarrow{}\infty$ and $\displaystyle \boldsymbol{H_t}\xrightarrow{a.s.}\boldsymbol{R}$ as $t\xrightarrow{}\infty$. Thus $\displaystyle \boldsymbol{H_t}\boldsymbol{G_t}^{-1}\xrightarrow{a.s.}\boldsymbol{R}$ as $t\xrightarrow{}\infty$. We also observe that $\boldsymbol{D_t}$ is bounded. This is because, $\displaystyle \|\boldsymbol{D_t}\|_{\ell_1}=\sum_{i=1}^m \|\boldsymbol{A_t^i}\|_{\ell_1}=m$. Thus $\Bar{\boldsymbol{D_t}}$ is also bounded. Hence, we can obtain the following convergence,
\begin{equation}
    \label{eq:L1convergence}
    \Bar{\boldsymbol{D_t}}(\boldsymbol{H_t}\boldsymbol{G_t}^{-1}-\boldsymbol{R})\xrightarrow{a.s.}0  \;\; \text{as}\;t\xrightarrow{}\infty.
\end{equation}
We now simply add the equations \eqref{eq:combined equation} and \eqref{eq:L1convergence} to obtain,
   \begin{equation}
       \label{eq:convergence of Dt}
       \boldsymbol{D_{t+1}}-\Bar{\boldsymbol{D_t}}\boldsymbol{R} \xrightarrow{a.s.}0 \;\; \text{as}\;t\xrightarrow{}\infty.
   \end{equation}
We know that $\boldsymbol{R}$ is a non-negative irreducible stochastic matrix (by the assumption made in \Cref{thm:deterministic interacting}) and $\{\boldsymbol{D_t}\}_{t\in\mathbb{N}}$ is a sequence of non-negative normalized vectors. Suppose, $\boldsymbol{\pi_{R}}$ is the dominant eigen-vector of $\boldsymbol{R}$. Hence, on applying \Cref{lem:gouet} to \eqref{eq:convergence of Dt}, we get,
\begin{equation}
    \boldsymbol{D_t}\xrightarrow{a.s.}\boldsymbol{\pi_R} \;\;\text{as}\;t\xrightarrow{}\infty.
\end{equation}
Let $\boldsymbol{\pi_1}$ denote the sub-vector formed by the first $N$ entries of $\boldsymbol{\pi_{R}}$, $\boldsymbol{\pi_2}$ denote the sub-vector formed by the next $N$ entries of $\boldsymbol{\pi_{R}}$ and so on. Thus, we have, $\boldsymbol{\pi}=(\boldsymbol{\pi_1},\boldsymbol{\pi_2},\cdots,\boldsymbol{\pi_m})$. We note that as $\displaystyle \boldsymbol{D_t}\xrightarrow{a.s.}\boldsymbol{\pi_{R}}$ as $t\xrightarrow{}\infty$, we have $\displaystyle \boldsymbol{A_t^i} \xrightarrow{a.s.}\boldsymbol{\pi_i}$ as $t\xrightarrow{}\infty$ for $i=1,2,\cdots,m$. From the equation, $\boldsymbol{\pi_{R}}\boldsymbol{R}=\boldsymbol{\pi_{R}}$, we obtain the following $m$ relations, $\boldsymbol{\pi_m}\boldsymbol{R}_1=\boldsymbol{\pi_1}$, $\boldsymbol{\pi_1}\boldsymbol{R}_2=\boldsymbol{\pi_2}$,$\cdots$,$\boldsymbol{\pi_{m-1}}\boldsymbol{R}_m=\boldsymbol{\pi_m}$. From these $m$ equations, we obtain that $\displaystyle \boldsymbol{\pi_1}\boldsymbol{R}_2\boldsymbol{R}_3\cdots\boldsymbol{R}_m\boldsymbol{R}_1=\boldsymbol{\pi_2}\boldsymbol{R}_3\cdots\boldsymbol{R}_m\boldsymbol{R}_1=\cdots=\boldsymbol{\pi_1}$. Thus, $\boldsymbol{\pi_1}$ is the dominant eigen-vector of the $N \times N$ matrix $\boldsymbol{R}_2\boldsymbol{R}_3\cdots\boldsymbol{R}_m\boldsymbol{R}_1$. It can be shown similarly that $\boldsymbol{\pi_i}$ is the dominant eigen-vector of the $N \times N$ matrix $\boldsymbol{R}_{i+1}\boldsymbol{R}_{i+2}\cdots\boldsymbol{R}_m\boldsymbol{R}_1\boldsymbol{R}_2\cdots\boldsymbol{R}_i$ for $i=1,2,\cdots,m$. Hence we can say that $\displaystyle \boldsymbol{\pi_i}=\boldsymbol{\pi_{{R}_{i+1}{R}_{i+2}\cdots{R}_m{R}_1{R}_2\cdots{R}_i}}$ for $i=1,2,\cdots,m$. This completes the proof.
\end{proof}
We note that one of the main assumptions of \textbf{Theorem 2} is that the $Nm \times Nm$ matrix, $\boldsymbol{R}$ is irreducible. This might be a difficult condition to verify in some situations. We thus provide below a necessary and sufficient condition to verify that $\boldsymbol{R}$ is an irreducible matrix.
\begin{lemma}\label{lem:iff condition for irreducibility}
$\boldsymbol{R}$ is an irreducible matrix iff $\boldsymbol{R}_{i+1}\boldsymbol{R}_{i+2}\cdots\boldsymbol{R}_m\boldsymbol{R}_1\cdots\boldsymbol{R}_i$ is an irreducible matrix for $i=1,2,...,m$.
\end{lemma}
\begin{proof}[Proof of \Cref{lem:iff condition for irreducibility}]
At first, let us prove the 'only if' part. Thus, we need to prove that $\boldsymbol{R}_{i+1}\boldsymbol{R}_{i+2}\cdots\boldsymbol{R}_m\boldsymbol{R}_1\cdots\boldsymbol{R}_i$ is an irreducible matrix $\forall \; i$ assuming that $\boldsymbol{R}$ is an irreducible matrix. Let us take a look at the powers of $\boldsymbol{R}$. $\boldsymbol{R}^2$ has the following form,
\begin{equation}
    \label{eq:2nd power of R}
    \boldsymbol{R}^2=\begin{pmatrix} 0&0&\boldsymbol{R}_2\boldsymbol{R}_3&0&\hdots&0 \\
    0&0&0&\boldsymbol{R}_3\boldsymbol{R}_4&\hdots&0\\
    \vdots&\vdots&\vdots&\vdots&\ddots&\vdots\\
    \boldsymbol{R}_m\boldsymbol{R}_1&0&0&0&\hdots&0\\
    0&\boldsymbol{R}_1\boldsymbol{R}_2&0&0&\hdots&0
    \end{pmatrix}.
\end{equation}
Continuing like this, $\boldsymbol{R}^m$ will have the following form,
\begin{equation}
    \label{eq:mth power of R}
    \boldsymbol{R}^m=\begin{pmatrix} \boldsymbol{R}_2\boldsymbol{R}_3...\boldsymbol{R}_1&0&0&\hdots&0 \\
    0&\boldsymbol{R}_3\boldsymbol{R}_4...\boldsymbol{R}_2&0&\hdots&0\\
    \vdots&\vdots&\vdots&\ddots&\vdots\\
    0&0&0&\hdots&\boldsymbol{R}_1\boldsymbol{R}_2...\boldsymbol{R}_m
    \end{pmatrix}.
\end{equation}
Now, since $\boldsymbol{R}$ is an irreducible matrix, this implies that for $1 \leq i,j \leq Nm$, $\exists$ a natural number $k$ such that $(\boldsymbol{R}^k)_{ij}>0$. Note that the diagonal blocks of $\boldsymbol{R}^k$ are non-zero only when $k=m, 2m, 3m,\cdots$. Moreover, we note that when $k=m$, the diagonal blocks are $\boldsymbol{R}_{i+1}\cdots\boldsymbol{R}_i$ where $i=1,2,\cdots,m$. Similarly, when $k=tm$, the diagonal blocks are $(\boldsymbol{R}_{i+1}\cdots\boldsymbol{R}_i)^t$ where $i=1,2,\cdots,m$. Hence, for $\boldsymbol{R}$ to be an irreducible matrix, $\boldsymbol{R}_{i+1}\boldsymbol{R}_{i+2}\cdots\boldsymbol{R}_m\boldsymbol{R}_1...\boldsymbol{R}_i$ must be an irreducible matrix for $i=1,2,\cdots,m$. This completes the proof of the only if (necessary condition) part.

Now let us prove the 'if' part. Here, we are given that $\boldsymbol{R}_{i+1}\boldsymbol{R}_{i+2}\cdots\boldsymbol{R}_m\boldsymbol{R}_1\cdots\boldsymbol{R}_i$ is an irreducible matrix for $i=1,2,\cdots,m$. We need to show that $\boldsymbol{R}$ is an irreducible matrix. Before proving this, we need to prove two claims.
\begin{claim}
\label{clm:UV}
Suppose, $\boldsymbol{U} \And \boldsymbol{V}$ are both $n \times n$ matrices, have non-negative entries and are balanced (each row sum is equal to 1). Suppose $\boldsymbol{V}$ is an irreducible matrix. Then, given $1 \leq i,j \leq n$, $\exists$ $k\in \mathbb{N}$ such that, $(\boldsymbol{U}\boldsymbol{V}^k)_{ij}>0$.
\end{claim}
\begin{proof}[Proof of \Cref{clm:UV}]
 We note that, $(\boldsymbol{U}\boldsymbol{V}^k)_{ij}=\boldsymbol{U}_{i*}\boldsymbol{V}^k_{*j}$. Since, sum of all entries in $\boldsymbol{U}_{i*}$ is 1, $\exists$ $1\leq s \leq n$ such that, $\boldsymbol{U}_{is}>0$. Also, since $\boldsymbol{V}$ is irreducible, $\exists$ $k \in N$ such that $\boldsymbol{V}^k_{sj}>0$. Thus, for this $k$ we have, $\displaystyle (\boldsymbol{U}\boldsymbol{V}^k)_{ij}=\boldsymbol{U}_{i*}\boldsymbol{V}^k_{*j} \geq \boldsymbol{U}_{is}\boldsymbol{V}^k_{sj} > 0$.
 \end{proof}
 \begin{claim}
 \label{clm:XY}
 Suppose $\boldsymbol{X},\boldsymbol{Y}$ are two $n \times n$ matrices with non-negative entries and are balanced (each row sum is equal to 1). Then $\boldsymbol{XY}$ satisfies the same properties as $\boldsymbol{X}$ i.e. it is also balanced (each row sum is $1$).
 \end{claim}
\Cref{clm:XY} can be proved simply by multiplying $X$ and $Y$. Now, let's come back to the proof of the sufficiency part. Note that any sub-matrix in $R^k$ has the following form: $\boldsymbol{R}_{i_1}\boldsymbol{R}_{i_2}\cdots\boldsymbol{R}_{i_p}(\boldsymbol{R}_{l+1}\cdots\boldsymbol{R}_m\cdots\boldsymbol{R}_l)^q$ for some $p,l,q \in \mathbb{N}$. By \Cref{clm:XY}, we know that $\boldsymbol{R}_{i_1}\boldsymbol{R}_{i_2}\cdots\boldsymbol{R}_{i_p}$ is balanced (each row sum is 1). Setting $\boldsymbol{U}=\boldsymbol{R}_{i_1}\boldsymbol{R}_{i_2}\cdots\boldsymbol{R}_{i_p}$ and $\boldsymbol{V}=\boldsymbol{R}_{l+1}\cdots\boldsymbol{R}_m\cdots\boldsymbol{R}_l$ (irreducible) and thereafter using \Cref{clm:UV}, we are done. This completes the proof of the if (sufficient condition) part and hence also the proof of \Cref{lem:iff condition for irreducibility}.
\end{proof}
We shall see some applications of \Cref{thm:deterministic interacting} to urn models in \Cref{sec:consequences}.

\section{Feedback Interacting Urn Model with Non-Deterministic Interaction}
\label{sec:Non-deterministic urn}
 Instead of working with a deterministic model as in \Cref{sec:deterministic cyclic}, here we will work with a non-deterministic interacting urn model. Like in the previous case, we have $m$ urns and each of them may contain balls of $N$ different colors. The replacement matrices of the $m$ urns are given by $\boldsymbol{R}_1, \boldsymbol{R}_2, \cdots,\boldsymbol{R}_m$, all of which are $N \times N$ matrices. Also, there is an $m \times m$ stochastic matrix, $\boldsymbol{P}$ given by,
 \begin{equation}
    \label{eq:P matrix}
    \boldsymbol{P}=\begin{pmatrix} p_{11} & p_{12} &\hdots&p_{1m} \\
    p_{21} & p_{22} &\hdots& p_{2m}\\
    \vdots&\vdots&\ddots&\vdots \\
    p_{m1} & p_{m2} &\hdots& p_{mm}
    \end{pmatrix}.
\end{equation}
It is assumed that all the urns have at least one ball in them initially. Now the following iterative process is performed:
\begin{enumerate}
    \item A ball is chosen from the first urn. The color of the ball is noted. Then, an urn is randomly chosen from the $m$ urns based on the probability vector $P_{1.}$ i.e. the first row of $P$. Hence the first urn is chosen with probability $p_{11}$, the second urn is chosen with probability $p_{12}$ and so on. Then, balls are added to the chosen urn according to the replacement matrix of that urn based on the color of the ball chosen from the first urn.
    \item A ball is chosen from the second urn. The color of the ball is noted. Then, an urn is randomly chosen from the $m$ urns based on the probability vector $P_{2.}$ i.e. the second row of $P$. Then, balls are added to the chosen urn according to the replacement matrix of that urn based on the color of the ball chosen from the second urn.
    \item $\cdots$
    \item A ball is chosen from the $m$-th urn. The color of the ball is noted. Then, an urn is randomly chosen from the $m$ urns based on the probability vector $P_{m.}$ i.e. the $m$-th row of $P$. Then, balls are added to the chosen urn according to the replacement matrix of that urn based on the color of the ball chosen from the $m$-th urn.
    \item Go to step(1).
\end{enumerate}

We have the following theorem for this set-up.
\begin{theorem}\label{thm:non deterministic urn}
As in \Cref{thm:deterministic interacting}, let $\boldsymbol{A_k^i}$ denote the proportion vector of urn-$i$ at the $k$-th stage. Let $\displaystyle \Tilde{\boldsymbol{R}}$ be a $Nm \times Nm$ matrix defined as follows, \begin{equation}
    \label{eq:R tilde}
    \Tilde{\boldsymbol{R}}=\begin{pmatrix}p_{11}\boldsymbol{R_1}&p_{12}\boldsymbol{R_2}&\hdots&p_{1m}\boldsymbol{R_m}\\
    p_{21}\boldsymbol{R_1}&p_{22}\boldsymbol{R_2}&\hdots&p_{2m}\boldsymbol{R_m}\\
    \vdots&\vdots&\ddots&\vdots\\
     p_{m1}\boldsymbol{R_1}&p_{m2}\boldsymbol{R_2}&\hdots&p_{mm}\boldsymbol{R_m}
    \end{pmatrix}.
\end{equation}
Let $\boldsymbol{\pi_{\Tilde{R}}}$ be the left eigen-vector of $\Tilde{\boldsymbol{R}}$ corresponding to the eigen-value $1$. Suppose $\boldsymbol{\pi_1}$ denotes the sub-vector formed by the first $N$ entries of $\boldsymbol{\pi_{\Tilde{R}}}$, $\boldsymbol{\pi_2}$ denotes the sub-vector formed by the next $N$ entries of $\boldsymbol{\pi_{\Tilde{R}}}$ and so on. We make the following two assumptions: \textit{(i)} $\boldsymbol{R}_i$ is a balanced (each row sum is 1) replacement matrix with all the entries non-negative for i=1,2,$\cdots$,m and \textit{(ii)} $\Tilde{\boldsymbol{R}}$ is an irreducible matrix. If we have the set-up described as above, then, we have the following result:
\begin{equation}
\nonumber
    \boldsymbol{A_p^i} \xrightarrow{a.s.} \boldsymbol{\pi_i} \;\;\text{as}\;\; p\xrightarrow{}\infty .
\end{equation}
\end{theorem}
\begin{proof}[Proof of \Cref{thm:non deterministic urn}]
We shall construct a $m^2N\times m^2N$ combined replacement matrix $\boldsymbol{R}$ from these $m$ matrices. Let $\boldsymbol{R}_{ij}$ denote the $ij$-th $Nm \times Nm$ block of $R$. The matrix $R$ is constructed in such a way that $[\boldsymbol{R}_{ij}]_{ji}=\boldsymbol{R}_j$ for $1\leq i,j \leq m$ and $0$ otherwise. Let us give an example to illustrate how the matrix $\boldsymbol{R}$ is constructed from the given $m$ replacement matrices. We give an example for the case when $m=3$ below,
\begin{equation}
    \boldsymbol{R}=\begin{pmatrix} \boldsymbol{R}_1&0&0&0&0&0&0&0&0\\
    0&0&0&\boldsymbol{R}_2&0&0&0&0&0\\
    0&0&0&0&0&0&\boldsymbol{R}_3&0&0\\
    0&\boldsymbol{R}_1&0&0&0&0&0&0&0\\
    0&0&0&0&\boldsymbol{R}_2&0&0&0&0\\
    0&0&0&0&0&0&0&\boldsymbol{R}_3&0\\
    0&0&\boldsymbol{R}_1&0&0&0&0&0&0\\
    0&0&0&0&0&\boldsymbol{R}_2&0&0&0\\
    0&0&0&0&0&0&0&0&\boldsymbol{R}_3
    \end{pmatrix}. \nonumber
\end{equation}
As in the earlier proof, we re-frame the given model by combining all the $m$ urns into a single urn having replacement matrix $\boldsymbol{R}$ as described above. Suppose the composition vector of the combined urn at the $t$-th stage is given by $\displaystyle \boldsymbol{C_t}=(\boldsymbol{C_t^{11}},\boldsymbol{C_t^{12}},\cdots,\boldsymbol{C_t^{1m}},\boldsymbol{C_t^{21}},\cdots,\boldsymbol{C_t^{2m}},\cdots,\boldsymbol{C_t^{mm}})$. Here, $\boldsymbol{C_t^{ij}}$ is the composition vector for the number of balls that have been added to the $i$-th urn due to color selection from the $j$-th urn till time $t$. So, $\boldsymbol{C_t^{ij}}$ is a $N$ dimensional vector for $1\leq i,j \leq m$. We also introduce the following notations: $\boldsymbol{C_t^1}=(\boldsymbol{C_t^{11}},\boldsymbol{C_t^{12}},\cdots,\boldsymbol{C_t^{1m}})$, $\boldsymbol{C_t^2}=(\boldsymbol{C_t^{21}},\boldsymbol{C_t^{22}},\cdots,\boldsymbol{C_t^{2m}})$ and so on. Thus $\boldsymbol{C_t^1},\boldsymbol{C_t^2},\cdots$ all are $mN$ dimensional row vectors. Let us denote by $\boldsymbol{\chi_n}$, the vector to indicate which color is chosen at the $n$th stage from the combined urn. Note that, $\boldsymbol{\chi_n}$ is a $m^2N$ dimensional row vector. Hence, $\boldsymbol{\chi_n}=\boldsymbol{e_i}$ if the $i$-th color is chosen at the $n$-th stage, where $i=1,2,3,\cdots,m^2N$ ($\{\boldsymbol{e_1}, \boldsymbol{e_2}, \boldsymbol{e_3},\cdots,\boldsymbol{e_{m^2N}}\}$ is the canonical basis of $\mathbb{R}^{m^2N}$). We need to keep in mind that in the first step, the ball is chosen from urn-1, in the second step, the ball is chosen from urn-2 and so on in a cyclic manner. Suppose, $\mathscr{F}_n$ is the sigma field generated by the set $\{\boldsymbol{\chi_i} \; for \; i=1,2,\cdots,n-1\}$. Therefore, the distribution of $\boldsymbol{\chi_n}$ is as described below. Let $S=\{(x,y)|x=1,2,\cdots,N \text{ and } y=0,1,2,\cdots,m-1\}$. For $n=mk+r$ where $r=1,2,\cdots,m$ we have,
\begin{equation}
\label{eq: chin distrn for non deterministic}
    \mathbb{P}(\boldsymbol{\chi_n}=\boldsymbol{e_{N(m(r-1)+t)+i}}|\mathscr{F}_n)=\begin{cases} &p_{r(t+1)}\sum_{j=0}^{m-1}\frac{\boldsymbol{C_{mk+r-1,Nj+i}}}{k+1} \; \text{for} \; (i,t)\in S ,\\ &0 \;\; \text{otherwise} .\end{cases}
\end{equation}
We shall now apply the Martingale version of Second Borel-Cantelli Lemma separately $m$ times on $\boldsymbol{\chi_{n,i}}$ depending on whether $n$ is equal to $mk+1, mk+2,\cdots, mk+m$ for some $k \in \mathbb{Z}^+ \cup \{0\}$. (Note that if $\boldsymbol{V_n}$ is a vector, then $\boldsymbol{V_{n,i}}$ denotes the $i$-th component of $\boldsymbol{V_n}$). We get the following $m$ equations using the distribution of $\boldsymbol{\chi_n}$ given in \eqref{eq: chin distrn for non deterministic}. For $n=mk+r$ where $r=1,2,\cdots,m$ we have,
\begin{equation}
\label{eq:bc on chin nondetr}
\frac{\sum_{k=0}^l \boldsymbol{\chi_{mk+r,Nt+i}}}{p_{r(t+1)}\sum_{k=0}^l\sum_{j=0}^{m-1} \frac{\boldsymbol{C_{mk+r-1,Nj+i}}}{k+1}} \xrightarrow{a.s.} 1 \; \text{on} \; \left\{ \sum_{k=0}^\infty \sum_{j=0}^{m-1}\frac{\boldsymbol{C_{mk+r-1,Nj+i}}}{k+1}=\infty \right\} \;\; \text{for} \; (i,t)\in S .
\end{equation}
We have, $\displaystyle \sum_{k=0}^\infty \sum_{j=0}^{m-1}{\boldsymbol{C_{mk+1,Nj+i}}}/({k+1}) \geq \sum_{k=0}^\infty {1}/({k+1})=\infty \; for \; i=1,2,\cdots,N $. Let $\boldsymbol{\chi_n^1}$ denote the sub-vector of $\boldsymbol{\chi_n}$ formed by the first $Nm$ entries, $\boldsymbol{\chi_n^2}$ denote the sub-vector of $\boldsymbol{\chi_n}$ formed by the next $Nm$ many entries and so on. From \eqref{eq:bc on chin nondetr}, we get the following set of $m$ equations in vector notation. For $n=mk+r$ where $r=1,2,\cdots,m$ we have,
\begin{equation}
    \sum_{k=0}^l \boldsymbol{\chi_{mk+r}^r}=\sum_{k=0}^l \frac{\boldsymbol{C_{mk+r-1}^r}}{k+1}.\paren{\boldsymbol{M^r} + o^r\paren{n}},
\end{equation}
where $\boldsymbol{M^r}$ is defined as,
\begin{equation}
    \label{eq:M^r}
    \boldsymbol{M^r}=\begin{pmatrix}
p_{r1}\boldsymbol{I_N} & p_{r2}\boldsymbol{I_N} &\hdots& p_{rm}\boldsymbol{I_N}\\
p_{r1}\boldsymbol{I_N} & p_{r2}\boldsymbol{I_N} &\hdots& p_{rm}\boldsymbol{I_N}\\
\vdots&\vdots&\ddots&\vdots\\
p_{r1}\boldsymbol{I_N} & p_{r2}\boldsymbol{I_N} &\hdots& p_{rm}\boldsymbol{I_N}\\
\end{pmatrix}.
\end{equation}
Note that $\boldsymbol{M^r}$ is a $Nm \times Nm$ matrix for $r=1,2,\cdots,m$. We know that for $n\geq1$ the urn composition evolves as in \eqref{eq:evolution of urn composition}. There are $m$ cases depending on whether $n$ is equal to $mp+1, mp+2,\cdots, mp+m$ for some $p \in \mathbb{Z}^+ \cup \{0\}$. We shall show the detailed calculation for only the last case i.e. when $n=mp$. The calculations for the other cases can be done quite similarly. Suppose, $n=mp$ for some $p \in \mathbb{Z}^+$. Then, we have,
\begin{equation}
    \label{eq:case m non deterministic}
    \begin{split}
    \sum_{j=1}^n\boldsymbol{\chi_j} &=\paren{\sum_{k=0}^{p-1}\boldsymbol{\chi_{mk+1}^1},\sum_{k=0}^{p-1}\boldsymbol{\chi_{mk+2}^2},\cdots,\sum_{k=1}^p\boldsymbol{\chi_{mk}^m}}\\ &=\paren{\sum_{k=0}^{p-1}\frac{\boldsymbol{C_{mk}^1}}{k+1}\paren{\boldsymbol{M^1}+o^1\paren{n}},\sum_{k=0}^{p-1}\frac{\boldsymbol{C_{mk+1}^2}}{k+1}\paren{\boldsymbol{M^2}+o^2\paren{n}},\cdots,\sum_{k=1}^p\frac{\boldsymbol{C_{mk-1}^m}}{k+1}\paren{\boldsymbol{M^m}+o^m\paren{n}}} \\ &=\paren{\sum_{k=0}^{p-1}\frac{\boldsymbol{C_{mk}^1}}{k+1},\sum_{k=0}^{p-1}\frac{\boldsymbol{C_{mk+1}^2}}{k+1},\cdots,\sum_{k=1}^p\frac{\boldsymbol{C_{mk-1}^m}}{k+1}}\paren{\Tilde{\boldsymbol{M}}+o\paren{n}},
    \end{split}
\end{equation}
where $o(n)$ is $Nm^2 \times Nm^2$ matrix whose diagonal submatrices are $o^1(n),o^2(n),\cdots,o^m(n)$ and all other entries are $0$. Similarly $\Tilde{\boldsymbol{M}}$ is $Nm^2 \times Nm^2$ matrix whose diagonal submatrices are $\boldsymbol{M^1},\boldsymbol{M^2},\cdots,\boldsymbol{M^m}$ and all other entries are $0$. Therefore, from equation \eqref{eq:evolution of urn composition}, we have,
\begin{equation}
      \boldsymbol{C_{mp}}-\boldsymbol{C_0}=\paren{\sum_{k=0}^{p-1}\frac{\boldsymbol{C_{mk}^1}}{k+1},\sum_{k=0}^{p-1}\frac{\boldsymbol{C_{mk+1}^2}}{k+1},\cdots,\sum_{k=1}^p\frac{\boldsymbol{C_{mk-1}^m}}{k+1}}\Tilde{\boldsymbol{M}}\boldsymbol{R} .
\end{equation}
We note that, in the matrix product, $\Tilde{\boldsymbol{M}}\boldsymbol{R}$, all the $\boldsymbol{R}_{ij}$'s will be replaced by $\boldsymbol{M^i}\boldsymbol{R}_{ij}$. In, $\boldsymbol{M^i}\boldsymbol{R}_{ij}$ all the entries in the $i$-th column are $p_{ij}\boldsymbol{R}_j$ and all other columns are $0$. Thus the matrix product $\boldsymbol{M^i}\boldsymbol{R}_{ij}$ looks like the following,
\begin{equation}
    \boldsymbol{M^i}\boldsymbol{R}_{ij}=\begin{pmatrix}
0&\hdots&p_{ij}\boldsymbol{R}_j&\hdots&0\\
0&\hdots&p_{ij}\boldsymbol{R}_j&\hdots&0\\
\vdots&\ddots&\vdots&\ddots&\vdots\\
0&\hdots&p_{ij}\boldsymbol{R}_j&\hdots&0
\end{pmatrix}.
\end{equation}
Also, let us see what happens when we multiply $\boldsymbol{C_n^1}$ with $\boldsymbol{M^i}\boldsymbol{R}_{ij}$. We shall have, $\displaystyle \boldsymbol{C_n^1M^iR}_{ij}=(0,0,\cdots,p_{ij}\boldsymbol{C_n^{1t}},0,0,\cdots,0)\boldsymbol{R}_j$. Note that all the entries of $\boldsymbol{C_n^1M^iR}_{ij}$ are $0$ except the $i$-th entry, which is $p_{ij}\boldsymbol{C_n^{1t}R}_j$. Here, $\boldsymbol{C_n^{1t}}=\boldsymbol{C_n^{11}}+\boldsymbol{C_n^{12}}+\cdots+\boldsymbol{C_n^{1m}}$. Returning back to our proof, we have the following equation,
\begin{equation}
    \boldsymbol{C_{mp}^1}-\boldsymbol{C_0^1}=\paren{p_{11}\sum_{k=0}^{p-1}\frac{\boldsymbol{C_{mk}^{1t}}}{k+1}\boldsymbol{R_1},p_{21}\sum_{k=0}^{p-1}\frac{\boldsymbol{C_{mk+1}^{2t}}}{k+1}\boldsymbol{R_1},\cdots,p_{m1}\sum_{k=1}^p\frac{\boldsymbol{C_{mk-1}^{mt}}}{k+1}\boldsymbol{R_1}},
\end{equation}
which implies,
\begin{equation}
\label{eq:composition casem nondeter}
    \frac{\boldsymbol{C_{mp}^{1t}}-\boldsymbol{C_0^{1t}}}{p+1}=\frac{1}{p+1}\paren{p_{11}\sum_{k=0}^{p-1}\frac{\boldsymbol{C_{mk}^{1t}}}{k+1}+p_{21}\sum_{k=0}^{p-1}\frac{\boldsymbol{C_{mk+1}^{2t}}}{k+1}+\cdots+p_{m1}\sum_{k=1}^p\frac{\boldsymbol{C_{mk-1}^{mt}}}{k+1}}R_1.
\end{equation}
We remember from \Cref{thm:deterministic interacting} that the proportion vector of the urn-$i$ at the $k$-th stage is defined as, $\displaystyle \boldsymbol{A_k^i}=({\boldsymbol{C_{m(k-1)+i-1}}^{it}})/{k} $. Using these notations, equation \eqref{eq:composition casem nondeter} boils down to,
\begin{equation}
    \boldsymbol{A_{p+1}^1}-(p_{11}\overline{\boldsymbol{A_p^1}}+p_{21}\overline{\boldsymbol{A_p^2}}+\cdots+p_{m1}\overline{\boldsymbol{A_p^m}})\boldsymbol{R_1} \xrightarrow{a.s.}0 \;\; \text{as} \; p\xrightarrow{} \infty.
\end{equation}
Similarly, we get a total of $m$ almost sure convergence equations as $p\xrightarrow{} \infty$.
\begin{equation}
    \begin{split}
        &\boldsymbol{A_{p+1}^1}-(p_{11}\overline{\boldsymbol{A_p^1}}+p_{21}\overline{\boldsymbol{A_p^2}}+\cdots+p_{m1}\overline{\boldsymbol{A_p^m}})\boldsymbol{R_1} \xrightarrow{a.s.}0, \\
    &\boldsymbol{A_{p+1}^2}-(p_{12}\overline{\boldsymbol{A_{p+1}^1}}+p_{22}\overline{\boldsymbol{A_p^2}}+\cdots+p_{m2}\overline{\boldsymbol{A_p^m}})\boldsymbol{R_2} \xrightarrow{a.s.}0, \\
    &...\\
    &\boldsymbol{A_{p+1}^m}-(p_{1m}\overline{\boldsymbol{A_{p+1}^1}}+p_{2m}\overline{\boldsymbol{A_{p+1}^2}}+\cdots+p_{(m-1)m}\overline{\boldsymbol{A_{p+1}^{(m-1)}}}+p_{mm}\overline{\boldsymbol{A_p^m}})\boldsymbol{R_m} \xrightarrow{a.s.}0 .
    \end{split}
\end{equation}
We can rewrite these $m$ equations as follows,
\begin{equation}
     \label{eq:m equations non deter}
     \begin{split}
         &\boldsymbol{A_{p+1}^1}-(p_{11}\overline{\boldsymbol{A_p^1}}+p_{21}\overline{\boldsymbol{A_p^2}}+\cdots+p_{m1}\overline{\boldsymbol{A_p^m}})\boldsymbol{R_1} \xrightarrow{a.s.}0 ,\\
&-\frac{p_{12}\boldsymbol{R_2}}{p+1}\boldsymbol{A_{p+1}^1}+\boldsymbol{A_{p+1}^2}-(p_{12}\frac{p}{p+1}\overline{\boldsymbol{A_p^1}}+p_{22}\overline{\boldsymbol{A_p^2}}+\cdots+p_{m2}\overline{\boldsymbol{A_p^m}})\boldsymbol{R_2} \xrightarrow{a.s.}0 ,\\
&-\frac{p_{13}\boldsymbol{R_3}}{p+1}\boldsymbol{A_{p+1}^1}-\frac{p_{23}\boldsymbol{R_3}}{p+1}\boldsymbol{A_{p+1}^2}+\boldsymbol{A_{p+1}^3}-(p_{13}\frac{p}{p+1}\overline{\boldsymbol{A_p^1}}+p_{23}\frac{p}{p+1}\overline{\boldsymbol{A_p^2}}+p_{33}\overline{\boldsymbol{A_p^3}}\cdots+p_{m3}\overline{\boldsymbol{A_p^m}})\boldsymbol{R_3} \xrightarrow{a.s.}0 ,\\
     \end{split}
\end{equation}
and so on. Now, we shall combine all the $m$ equations into a single equation. Let $\displaystyle \boldsymbol{D_p}=(\boldsymbol{A_p^1},\boldsymbol{A_p^2},\cdots,\boldsymbol{A_p^m})$. So, on combining all the $m$ equations into a single equation, we get,
\begin{equation}
    \boldsymbol{D_{p+1}G_p}-\overline{\boldsymbol{D_p}}\boldsymbol{H_p} \xrightarrow{a.s.}0\;\;\text{as}\;\;p\xrightarrow{}\infty,
\end{equation}
which implies that,
\begin{equation}
\label{eq:Dp equation nondetr}
    \boldsymbol{D_{p+1}}-\overline{\boldsymbol{D_p}}\boldsymbol{H_pG_p}^{-1}\xrightarrow{a.s.}0\;\;\text{as}\;\;p\xrightarrow{}\infty,
\end{equation}
where,
\begin{equation}
    \label{eq:Gp}
   \boldsymbol{G_p}=\begin{pmatrix}\boldsymbol{I_N}&-\frac{p_{12}\boldsymbol{R_2}}{p+1}&-\frac{p_{13}\boldsymbol{R_3}}{p+1}&\hdots&-\frac{p_{1m}\boldsymbol{R_m}}{p+1}\\
        0&\boldsymbol{I_N}&-\frac{p_{23}\boldsymbol{R_3}}{p+1}&\hdots&-\frac{p_{2m}\boldsymbol{R_m}}{p+1}\\
        \vdots&\vdots&\vdots&\ddots&\vdots\\
        0&0&0&\hdots&\boldsymbol{I_N}
        \end{pmatrix},
\end{equation}
and,
\begin{equation}
\label{eq:Hp}
    \boldsymbol{H_p}=\begin{pmatrix}p_{11}\boldsymbol{R_1}&p_{12}\frac{p\boldsymbol{R_2}}{p+1}&p_{13}\frac{p\boldsymbol{R_3}}{p+1}&\hdots&p_{1m}\frac{p\boldsymbol{R_m}}{p+1}\\
        p_{21}\boldsymbol{R_1}&p_{22}\boldsymbol{R_2}&p_{23}\frac{p\boldsymbol{R_3}}{p+1}&\hdots&p_{2m}\frac{p\boldsymbol{R_m}}{p+1}\\
        \vdots&\vdots&\vdots&\ddots&\vdots\\
        p_{m1}\boldsymbol{R_1}&p_{m2}\boldsymbol{R_2}&p_{m3}\boldsymbol{R_3}&\hdots&p_{mm}\boldsymbol{R_m}
        \end{pmatrix}.
\end{equation}
We observe that, $\displaystyle \boldsymbol{G_p}\xrightarrow{a.s.}\boldsymbol{I_{Nm}}$ as $p\xrightarrow{}\infty$ and $\displaystyle \boldsymbol{H_p}\xrightarrow{a.s.} \Tilde{\boldsymbol{R}}$ as $p\xrightarrow{}\infty$. Thus $\displaystyle \boldsymbol{H_pG_p}^{-1}\xrightarrow{a.s.} \Tilde{\boldsymbol{R}}$ as $p\xrightarrow{}\infty$. We also observe that $\boldsymbol{D_p}$ is bounded. This is because, $\displaystyle \|\boldsymbol{D_p}\|_{\ell_1}=\sum_{i=1}^m \|\boldsymbol{A_p^i}\|_{\ell_1}=m$. Thus $\Bar{\boldsymbol{D_p}}$ is also bounded. Hence, we can obtain the following equation,
\begin{equation}
    \label{eq:boundedness of Dp nondetr}
    \Bar{\boldsymbol{D_p}}(\boldsymbol{H_pG_p}^{-1}-\Tilde{\boldsymbol{R}})\xrightarrow{a.s.}0 \;\; \text{as}\;p\xrightarrow{}\infty.
\end{equation}
We now simply add the equations \eqref{eq:Dp equation nondetr} and \eqref{eq:boundedness of Dp nondetr} to obtain,
 \begin{equation}
       \label{eq:final equation Dp nondetr}
       \boldsymbol{D_{p+1}}-\Bar{\boldsymbol{D_p}}\Tilde{\boldsymbol{R}} \xrightarrow{a.s.}0 \;\; \text{as}\;p\xrightarrow{}\infty.
   \end{equation}
We know that $\Tilde{\boldsymbol{R}}$ is a non-negative irreducible stochastic matrix (by the assumption made in \Cref{thm:non deterministic urn}) and $\{\boldsymbol{D_p}\}_{p\in\mathbb{N}}$ is a sequence of non-negative normalized vectors. Suppose, $\pi_{\Tilde{R}}$ is the dominant eigen-vector of $\Tilde{\boldsymbol{R}}$. Hence, on applying \Cref{lem:gouet} to \eqref{eq:final equation Dp nondetr}, we get,
\begin{equation}
    \boldsymbol{D_p}\xrightarrow{a.s.}\boldsymbol{\pi_{\Tilde{{R}}}} \;\;\text{as}\;p\xrightarrow{}\infty.
\end{equation}
Hence we obtain that, $\displaystyle \boldsymbol{A_p^i}\xrightarrow{a.s.}\boldsymbol{\pi_i}$ as $p\xrightarrow{}\infty$ for $i=1,2,\cdots,m$. This completes the proof.
\end{proof}
We shall see a nice application of \Cref{thm:non deterministic urn} to an urn model in \Cref{sec:consequences}.

\section{Some Interesting Consequences}
\label{sec:consequences}
In this section, we shall study some interesting consequences of \Cref{thm:deterministic interacting} and \Cref{thm:non deterministic urn}. In \Cref{subsec:2urn} and \Cref{subsec:3urn}, we shall see the application of \Cref{thm:deterministic interacting} to a special deterministic $2$-urn and $3$-urn model respectively. In \Cref{subsec:nurn}, we shall see the application of \Cref{thm:non deterministic urn} to a special non-deterministic $n$-urn model for any given natural number, $n \geq 2$.

\subsection{The $2$-Urn-$2$-Color Model}
\label{subsec:2urn}
We consider a 2-urn-2-color model. In this model, we are given two urns, both of which contain balls of two colors say red (color 1) and black (color 2). The replacement matrices of the urns are $\boldsymbol{R_1}$ and $\boldsymbol{R_2}$. Both $\boldsymbol{R_1}$ and $\boldsymbol{R_2}$ are $2 \times 2$ non-random stochastic replacement matrices. The results obviously extend to non-random replacement matrices with constant row sums by obvious rescaling. $\boldsymbol{R_1}$ and $\boldsymbol{R_2}$ are defined as follows,
\begin{equation}
\label{eq:R1,R2}
\boldsymbol{R_1}=\begin{pmatrix} \alpha & 1-\alpha \\ 0 & 1 \end{pmatrix}, \;\;\;\;\; \boldsymbol{R_2}=\begin{pmatrix} 1 & 0 \\ \beta & 1-\beta  \end{pmatrix}.
\end{equation}
Note that, if $\boldsymbol{R_1}$ and $\boldsymbol{R_2}$ are reducible but not the identity matrix, then after possibly interchanging the names of the colors, they can be converted into upper triangular matrix (respectively lower triangular) as described above. It is assumed that both the urns have at least one ball in them initially and the iterative process mentioned in \Cref{sec:deterministic cyclic} is performed. As mentioned already in \Cref{sec:introduction}, we say that a particular urn \textit{flips} if the composition of the dominant color of the urn becomes less than the average composition of the remaining colors. Under this set-up, we have the following corollary.
\begin{corollary}\label{cor:2urn result}
Under the set-up described above, we have the following:
\begin{enumerate}
    \item If $\boldsymbol{X_t}$ and $\boldsymbol{Y_t}$ denote the proportion vectors of colors of first and second urn respectively then we have,
    \begin{equation}
    \begin{split}
      &\boldsymbol{X_t} \xrightarrow{a.s.} \paren{\frac{\alpha-\alpha\beta}{1-\alpha\beta}, \frac{1-\alpha}{1-\alpha\beta}},\\
   &\boldsymbol{Y_t} \xrightarrow{a.s.} \paren{\frac{1-\beta}{1-\alpha\beta}, \frac{\beta-\alpha\beta}{1-\alpha\beta}} .
    \end{split}
    \end{equation}
    \item In the limiting case, either none of the urns flip or exactly one of them flips. However, both of the urns can not flip simultaneously in the limiting case.
\end{enumerate}
\end{corollary}

\begin{proof}[Proof of \Cref{cor:2urn result}]
It can be easily shown that $\boldsymbol{R_1R_2}$ and $\boldsymbol{R_2R_1}$ both are irreducible matrices. Hence, from \Cref{lem:iff condition for irreducibility}, we obtain that $\boldsymbol{R}$ (as mentioned in \Cref{thm:deterministic interacting}) is an irreducible matrix. Thus all the assumptions of \Cref{thm:deterministic interacting} are satisfied. Calculating the values of $\boldsymbol{\pi_{R_1R_2}}$ and $\boldsymbol{\pi_{R_2R_1}}$ and applying \Cref{thm:deterministic interacting} proves the first part of the corollary.

For the second part, note that urn-1 flips if composition of the first color becomes more than the composition of the second color in the limiting case. Similarly, urn-2 flips if composition of the second color becomes more than the composition of the first color in the limiting case. We define two sets, $S_1 \And S_2$, where $S_1=\{(\alpha,\beta)|0<\alpha,\beta<1 \And 2\alpha > 1+\alpha\beta\}$ and $S_2=\{(\alpha,\beta)|0<\alpha,\beta<1 \And2\beta > 1+\alpha\beta\}$. Thus, we obtain that, urn-$1$ flips iff $(\alpha,\beta)\in S_1$ and urn-$2$ flips iff $(\alpha,\beta)\in S_2$. We will now show that both the urns can't flip simultaneously in the limiting case. Suppose, both the urns do flip simultaneously in the limiting case. This implies that $\beta>1/(2-\alpha)$ which in turn implies that $\beta>({2-\beta})/({3-2\beta})$. On simplifying this, we obtain, $(\beta -1)^2 < 0$, which is clearly a contradiction. This proves the second part of the corollary.
\end{proof}

\begin{figure}[htp]
    \centering
    \includegraphics[width=8 cm]{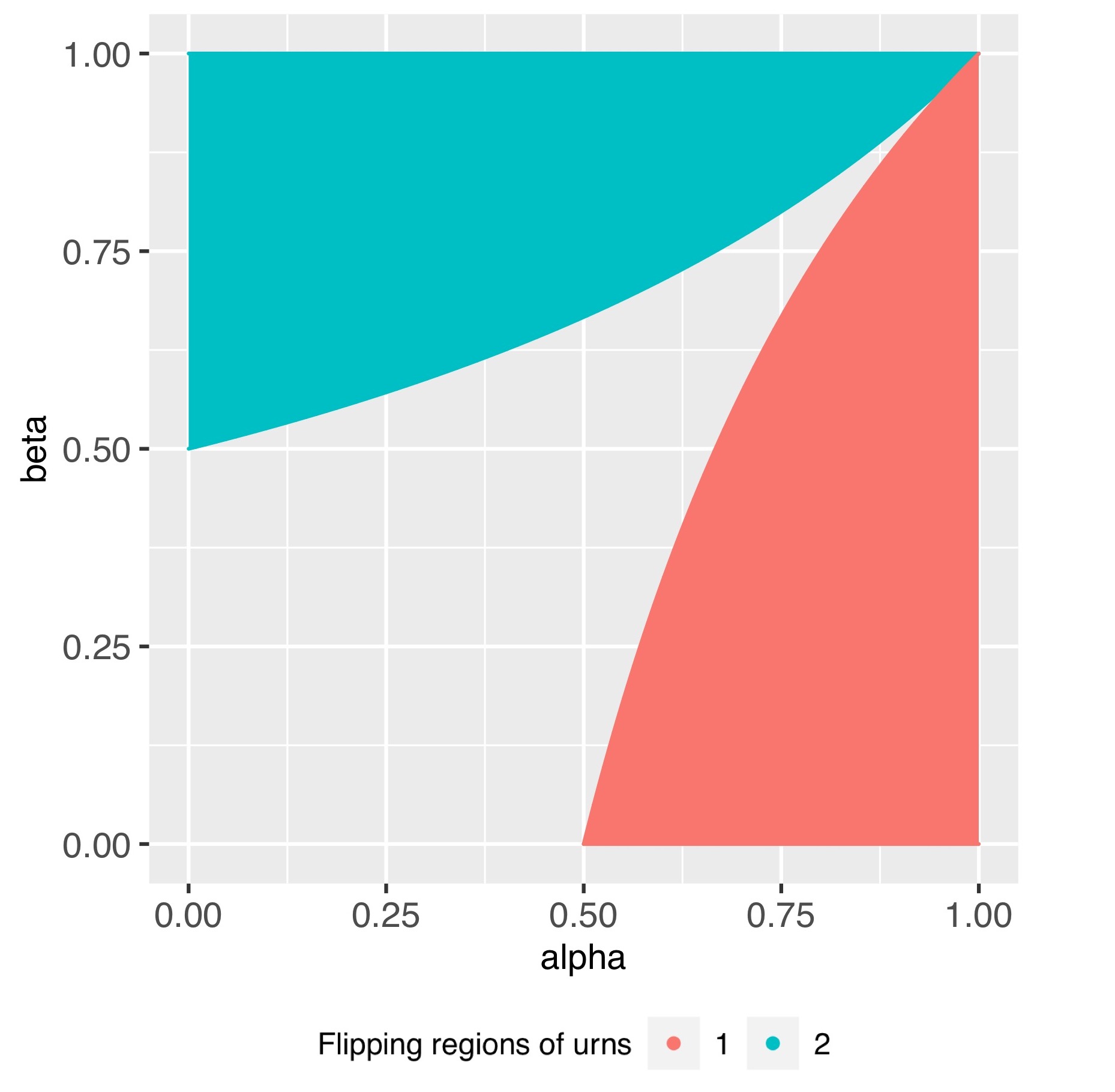}
    \caption{$2$-urn-$2$-color model.}
    \label{fig:2 urn}
\end{figure}

For better understanding, we are attaching the graph depicting the intersection of the two curves in \Cref{fig:2 urn}. Here, the wire-frame is that of a unit square in $2D$. The two curves are shown in two different colors. In the limiting case, the two urns flip respectively in the regions exterior to the two curves. We observe that the two urns can not flip simultaneously in the limiting case.

\subsection{The $3$-Urn-$3$-Color Model}
\label{subsec:3urn}
We consider a $3$-urn-$3$-color model. In other words, we have $m=N=3$. The replacement matrices of the urns, $\boldsymbol{R_1},\boldsymbol{R_2},\boldsymbol{R_3}$ are non-random stochastic matrices. Suppose, $\boldsymbol{R_1},\boldsymbol{R_2},\boldsymbol{R_3}$ are given by,
\begin{equation}
    \label{eq:R1,R2,R3}
    \boldsymbol{R_1}=\begin{pmatrix} \frac{\alpha}{2}&\frac{\alpha}{2} & 1-\alpha\\ \frac{\alpha}{2}& \frac{\alpha}{2}& 1-\alpha \\ 0&0&1
\end{pmatrix}, \; \boldsymbol{R_2}=\begin{pmatrix} \frac{\beta}{2}&1-\beta&\frac{\beta}{2}\\0&1&0\\ \frac{\beta}{2}& 1-\beta&\frac{\beta}{2}&
 \end{pmatrix}, \;\boldsymbol{R_3}=\begin{pmatrix} 1&0&0\\ 1-\gamma&\frac{\gamma}{2}&\frac{\gamma}{2}\\1-\gamma&\frac{\gamma}{2}&\frac{\gamma}{2} \end{pmatrix},
\end{equation}
where, $0<\alpha,\beta,\gamma<1$. It is assumed that all the urns have at least one ball in them initially and the iterative process mentioned in \Cref{sec:deterministic cyclic} is performed. Under this set-up, we have the following corollary.
\begin{corollary}\label{cor:3-urn}
Under the set-up described above, no two urns flip simultaneously in the limiting case.
\end{corollary}
\begin{proof}[Proof of \Cref{cor:3-urn}]
It can be easily shown that $\boldsymbol{R_1R_2R_3}$, $\boldsymbol{R_2R_3R_1}$ and $\boldsymbol{R_3R_1R_2}$ are all irreducible matrices. Hence, from \Cref{lem:iff condition for irreducibility}, we obtain that $\boldsymbol{R}$ (as mentioned in \Cref{thm:deterministic interacting}) is an irreducible matrix. Thus all the assumptions of \Cref{thm:deterministic interacting} are satisfied. Calculating the values of $\boldsymbol{\pi_{R_1R_2R_3}}$, $\boldsymbol{\pi_{R_2R_3R_1}}$, $\boldsymbol{\pi_{R_3R_1R_2}}$ and applying \Cref{thm:deterministic interacting} we obtain for $t\xrightarrow{}\infty$,
\begin{equation}
    \begin{split}
       & \boldsymbol{A_t^1} \xrightarrow{a.s.} \paren{ \frac{\alpha-\frac{\alpha\gamma}{2}+ \frac{\alpha\beta\gamma}{4}}{2+ \frac{\alpha\beta\gamma}{4}} ,
      \frac{\alpha-\frac{\alpha\gamma}{2}+ \frac{\alpha\beta\gamma}{4}}{2+ \frac{\alpha\beta\gamma}{4}} ,
     \frac{2-2\alpha + \alpha\gamma - \frac{\alpha\beta\gamma}{4}}{2+ \frac{\alpha\beta\gamma}{4}}} ,\\
     & \boldsymbol{A_t^2}\xrightarrow{a.s.}\paren{ \frac{\beta-\frac{\alpha\beta}{2}+ \frac{\alpha\beta\gamma}{4}}{2+ \frac{\alpha\beta\gamma}{4}},
      \frac{2-2\beta + \alpha\beta - \frac{\alpha\beta\gamma}{4}}{2+ \frac{\alpha\beta\gamma}{4}},
     \frac{\beta-\frac{\alpha\beta}{2}+ \frac{\alpha\beta\gamma}{4}}{2+ \frac{\alpha\beta\gamma}{4}}} ,\\
     &\boldsymbol{A_t^3}\xrightarrow{a.s.} \paren{\frac{2-2\gamma + \beta\gamma - \frac{\alpha\beta\gamma}{4}}{2+ \frac{\alpha\beta\gamma}{4}},      \frac{\gamma-\frac{\beta\gamma}{2}+ \frac{\alpha\beta\gamma}{4}}{2+ \frac{\alpha\beta\gamma}{4}} ,\frac{\gamma-\frac{\beta\gamma}{2}+ \frac{\alpha\beta\gamma}{4}}{2+ \frac{\alpha\beta\gamma}{4}} } .
    \end{split}
\end{equation}
Thus we obtain the following,
\begin{enumerate}
         \item $1$-st urn flips if $6\alpha-3\alpha\gamma+ \alpha\beta\gamma >4$.
         \item $2$-nd urn flips if $6\beta-3\alpha\beta+ \alpha\beta\gamma >4$.
         \item $3$-rd urn flips if $6\gamma-3\beta\gamma+ \alpha\beta\gamma >4$.
\end{enumerate}
We need to show that no two urns can flip simultaneously in the limiting case. Suppose, the $1$-st and $2$-nd urn flip simultaneously in the limiting case. This implies that the equations $6\alpha-3\alpha\gamma+ \alpha\beta\gamma >4$ and $6\beta-3\alpha\beta+ \alpha\beta\gamma >4$ hold simultaneously. From these two equations, we obtain,
\begin{equation}
    \beta > \frac{4}{6-\frac{4}{6-3\gamma +\beta\gamma}(3 -\gamma)},
\end{equation}
which implies,
\begin{equation}
    (\gamma)\beta^2+(4-3\gamma)\beta+(2\gamma-4)>0.
\end{equation}
This in turn implies that, $\beta<(2-({4}/{\gamma}))$ or $\beta>1$. As $0< \gamma <1$, we have, $\beta<-2$ or $\beta>1$. This is clearly a contradiction as we know that $0< \beta <1$. This completes the proof of the corollary.
\end{proof}

\begin{figure}[htp]
    \centering
    \includegraphics[width=8 cm]{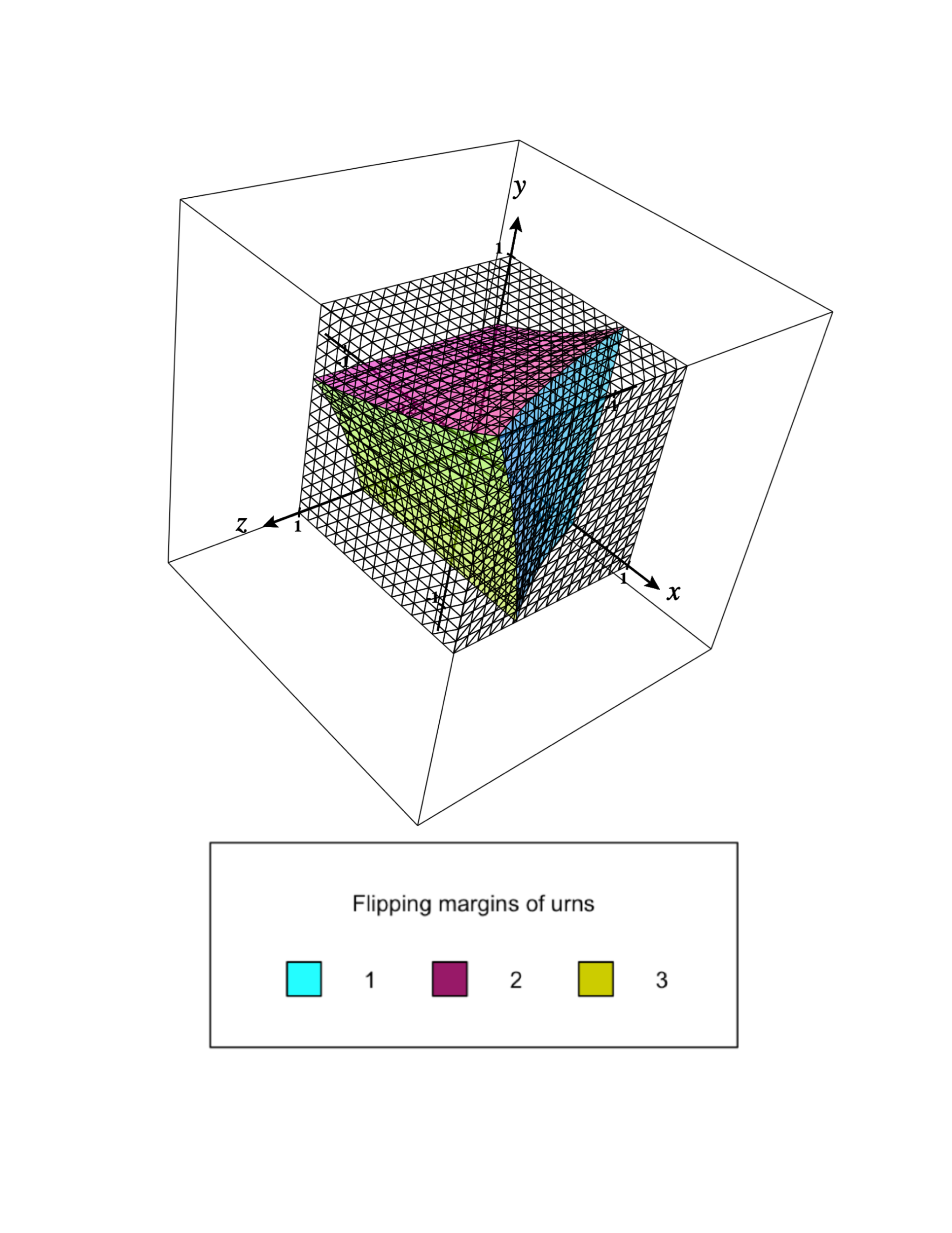}
    \caption{$3$-urn-$3$-color model.}
    \label{fig:3 urn}
\end{figure}

For better understanding, we are attaching the graph depicting the intersection of the three curves in \Cref{fig:3 urn}. In \Cref{fig:3 urn}, the wire-frame is that of a unit square in $3D$. The three curves are shown in three different colors. In the limiting case, the three urns flip respectively in the regions exterior to the three curves. We observe that no two urns can flip simultaneously in the limiting case.

\begin{figure}[htp]
    \centering
    \includegraphics[width=8 cm]{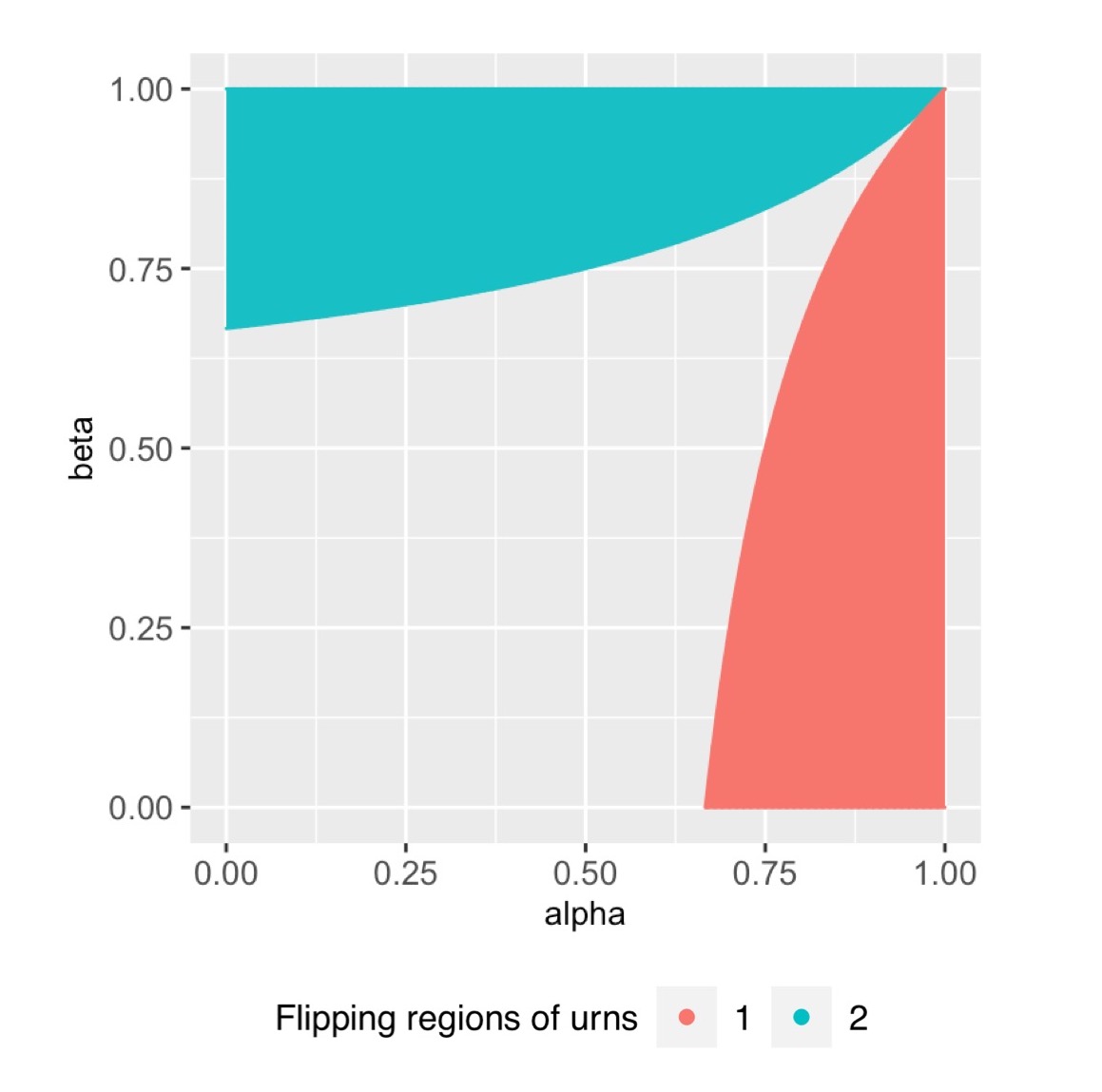}
    \caption{$n=2$ case of $n$-Urn-$n$-Color Model.}
    \label{fig:2 urn n}
\end{figure}

\subsection{The $n$-Urn-$n$-Color Model }
\label{subsec:nurn}
We consider a $n$-urn-$n$-color model i.e. we have $m=N=n$. The replacement matrices of the urns are non-random stochastic matrices. The replacement matrices of the $n$ urns are defined as follows,
\begin{equation}
    \label{eq:n urns non detr}
    \boldsymbol{R_1}=\begin{pmatrix} 1 & 0 & 0 &\hdots&0 \\
    1-a_1 &a_1&0&\hdots&0\\
    1-a_1&0&a_1&\hdots&0\\
    \vdots&\vdots&\vdots&\ddots&\vdots\\
    1-a_1&0&0&\hdots&a_1
    \end{pmatrix}, \; \; \;
     \boldsymbol{R_2}=\begin{pmatrix} a_2 & 1-a_2 & 0 &\hdots&0 \\
    0 &1&0&\hdots&0\\
    0&1-a_2&a_2&\hdots&0\\
    \vdots&\vdots&\vdots&\ddots&\vdots\\
    0&1-a_2&0&\hdots&a_2
    \end{pmatrix},
\end{equation}
and so on, where, $0<a_1,a_2,\cdots,a_n<1$. Note that, the replacement matrices have been defined in such a way that color-1 is the dominant color of the first urn, color-2 is the dominant color of the second urn and so on. In the stochastic matrix $\boldsymbol{P}$ (as defined in \Cref{sec:Non-deterministic urn}), all the entries $p_{ij}$ are taken to be $\frac{1}{n}$ $\forall$ $i\And j$. It is assumed that all the urns have at least one ball in them initially and the iterative process mentioned in \Cref{sec:Non-deterministic urn} is performed. Under this set-up we have the following corollary.
\begin{corollary}\label{cor:n urn}
In the limiting case, depending on the values of $a_1, a_2,\cdots,a_n$,  $k$ urns may flip simultaneously where $k=0,1,2,\cdots,n-1$. However, all the $n$ urns can never flip simultaneously in the limiting case.
\end{corollary}
\begin{proof}[Proof of \Cref{cor:n urn}]
We shall apply \Cref{thm:non deterministic urn} here. We know that, $\displaystyle \boldsymbol{A_p^i}\xrightarrow{a.s.}\boldsymbol{\pi_i}$ as $p\xrightarrow{}\infty$ for $i=1,2,\cdots,n$. We note that, in our case we have,
\begin{equation}
    \Tilde{\boldsymbol{R}}=\begin{pmatrix}
\frac{\boldsymbol{R_1}}{n}&\frac{\boldsymbol{R_2}}{n}&\hdots&\frac{\boldsymbol{R_n}}{n}\\
\frac{\boldsymbol{R_1}}{n}&\frac{\boldsymbol{R_2}}{n}&\hdots&\frac{\boldsymbol{R_n}}{n}\\
\vdots&\vdots&\ddots&\vdots\\
\frac{\boldsymbol{R_1}}{n}&\frac{\boldsymbol{R_2}}{n}&\hdots&\frac{\boldsymbol{R_n}}{n}
\end{pmatrix}.
\end{equation}
Let $\displaystyle \boldsymbol{\pi_{\Tilde{R}}}= ( \boldsymbol{\pi_1},\boldsymbol{\pi_2},\cdots,\boldsymbol{\pi_n})$. So we shall have,
\begin{equation}
    (\boldsymbol{\pi_1},\boldsymbol{\pi_2},\cdots,\boldsymbol{\pi_n}) \begin{pmatrix}
\frac{\boldsymbol{R_1}}{n}&\frac{\boldsymbol{R_2}}{n}&\hdots&\frac{\boldsymbol{R_n}}{n}\\
\frac{\boldsymbol{R_1}}{n}&\frac{\boldsymbol{R_2}}{n}&\hdots&\frac{\boldsymbol{R_n}}{n}\\
\vdots&\vdots&\ddots&\vdots\\
\frac{\boldsymbol{R_1}}{n}&\frac{\boldsymbol{R_2}}{n}&\hdots&\frac{\boldsymbol{R_n}}{n}
\end{pmatrix}
= (\boldsymbol{\pi_1},\boldsymbol{\pi_2},\cdots,\boldsymbol{\pi_n}).
\end{equation}
Thus, we shall get the following $n$ equations. For $i=1,2,\cdots,n$, we have,
\begin{equation}
    (\boldsymbol{\pi_1}+\boldsymbol{\pi_2}+\cdots+\boldsymbol{\pi_n})\boldsymbol{R_i}=n\boldsymbol{\pi_i}.
\end{equation}
Let $\boldsymbol{\pi_s}=\boldsymbol{\pi_1}+\boldsymbol{\pi_2}+\cdots\boldsymbol{\pi_n}$. Adding all the $n$ equations, we shall get, $\boldsymbol{\pi_s}\Bar{\boldsymbol{R}}=\boldsymbol{\pi_s}$, where $\Bar{\boldsymbol{R}}=({\boldsymbol{R_1}+\cdots+\boldsymbol{R_n}})/{n}$. We note that $\boldsymbol{\pi_s}$ is left eigen-vector of $\Bar{\boldsymbol{R}}$ corresponding to eigen-value $1$. Once $\boldsymbol{\pi_s}$ is obtained, we can obtain all the $\boldsymbol{\pi_i}$'s from the equation, $ \boldsymbol{\pi_i}= \boldsymbol{\pi_s} \boldsymbol{R_i} $ for $i=1,2,\cdots,n$. Note that,
\begin{equation}
    \Bar{\boldsymbol{R}}=\frac{1}{n}\begin{pmatrix}1+a_1+...+a_n & 1-a_2 & 1-a_3 & \hdots & 1-a_n\\
1-a_1 & 1+a_1+...+a_n & 1-a_3 &\hdots& 1-a_n\\
\vdots&\vdots&\vdots&\ddots&\vdots\\
1-a_1 & 1-a_2 & 1-a_3 &\hdots& 1+a_1+...+a_n \end{pmatrix}.
\end{equation}
By simple calculation, it can be shown that,
\begin{equation}
    \boldsymbol{\pi_s} = \paren{ \frac{1-a_1}{n- (a_1+ ... +a_n)} , \frac{1-a_2}{n- (a_1+ ... +a_n)} ,\cdots, \frac{1-a_n}{n- (a_1+ ... +a_n)}}.
\end{equation}
Let us now consider the first urn. From the equation $ \pi_1= \pi_s \boldsymbol{R_1} $, we obtain,
\begin{equation}
     \boldsymbol{\pi_1}=\paren{ \frac{(1-a_1)(n-a_2-...-a_n)}{n- (a_1+ ... +a_n)} , \frac{(1-a_2)a_1}{n- (a_1+ ... +a_n)} , \cdots, \frac{(1-a_n)a_1}{n- (a_1+ ... +a_n)}}.
\end{equation}
Urn-$1$ flips if composition of the first color becomes less than the average composition of the other colors. Thus, urn-$1$ flips if $\displaystyle \boldsymbol{\pi_1^2}+\boldsymbol{\pi_1^3}+\cdots+\boldsymbol{\pi_1^n} > (n-1)\boldsymbol{\pi_1^1}$. Substituting the values of $\boldsymbol{\pi_1^i}$  for $i=1,2,\cdots,n$, we get that, urn-$1$ flips if,
\begin{equation}
    (n^2-1)a_1 +(n-1)(a_2+a_3+\cdots+a_n) - n(a_1a_2+ a_1a_3+ \cdots + a_1a_n) > n(n-1).
\end{equation}
Similarly, we can obtain the condition of flipping for the remaining $n-1$ urns. Thus, we shall get a total of $n$ equations, where each equation represents the condition of flipping for a particular urn. We want to show that all the urns can not flip simultaneously in the limiting case. We shall prove this by the method of contradiction. Suppose all the urns do flip simultaneously. Then all the $n$ equations will hold simultaneously. Adding all the $n$ equations we get the following,
\begin{equation}
    (n^2-1)\paren{\sum_{i=1}^n a_i} +(n-1)^2\paren{\sum_{i=1}^n a_i} - n\sum_{i \neq j}a_ia_j > n^2(n-1),
\end{equation}
which implies,
\begin{equation}
    \sum_{i<j}(1-a_i)(1-a_j) < 0.
\end{equation}
This is clearly a contradiction as we know that, $0< a_1,a_2,\cdots,a_n <1$. Thus, we conclude that all the $n$ urns can never flip simultaneously in the limiting case. This completes the proof of the corollary.
\end{proof}

\begin{figure}[htp]
    \centering
    \includegraphics[width=8 cm]{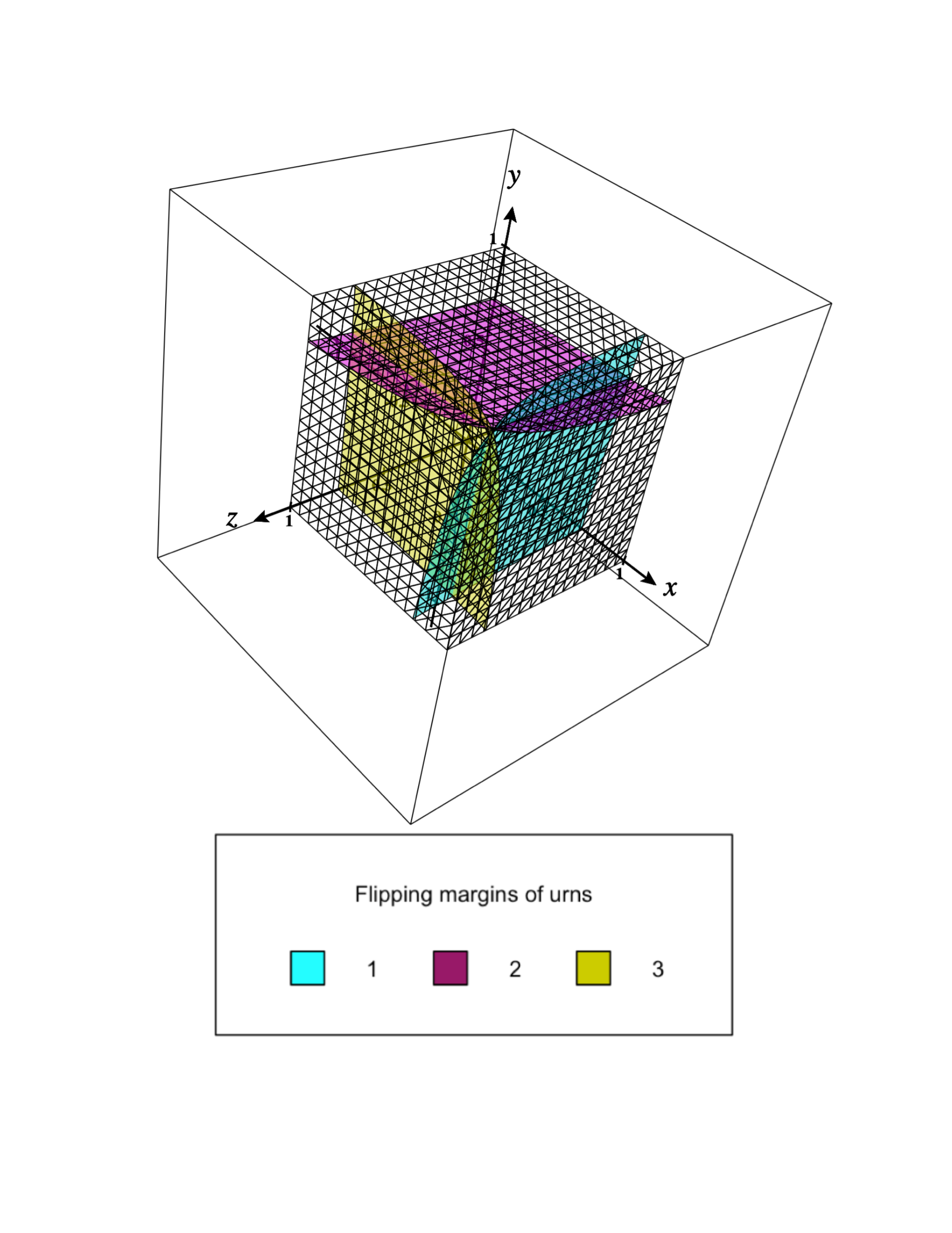}
    \caption{$n=3$ case of $n$-Urn-$n$-Color Model.}
    \label{fig:3 urn n}
\end{figure}

We shall give two graphical examples one for the case $n=2$ and the other for the case $n=3$. This shall help in better understanding and also serve as a justification of the correctness of \Cref{cor:n urn}. In both \Cref{fig:2 urn n} and \Cref{fig:3 urn n}, the wire-frame is that of a unit square and the curves are shown in different colors. As before, the urns flip in the regions exterior to the curves. We observe that in both \Cref{fig:2 urn n} and \Cref{fig:3 urn n}, all the urns can not flip simultaneously in the limiting case.




\section*{Acknowledgement}
The research of the first author was partly supported by MATRICS grant number MTR/2019/001448 from Science and Engineering Research Board, Govt.\ of India.

%
%
%
%
%
%
%
%
%
%
%

\bibliography{references}

\end{document}